\documentclass{amsart}
\usepackage{amssymb}
\newtheorem{theorem}{Theorem}[section]
\newtheorem{lemma}[theorem]{Lemma}
\newtheorem{proposition}[theorem]{Proposition}

\newtheorem{fact}[theorem]{Fact}

\theoremstyle{definition}

\theoremstyle{remark}

\numberwithin{equation}{section}

\renewcommand{\span}{\mathrm{span}}

\newcommand{\R}{\mathbb{R}}
\newcommand{\N}{\mathbb{N}}

\newcommand{\Q}{\mathbb{Q}}
\newcommand{\X}{\mathrm{X}}

\newcommand{\Y}{\mathrm{Y}}
\newcommand{\Z}{\mathrm{Z}}
\newcommand{\B}{\mathbf{B}}
\newcommand{\I}{\mathbf{I}}

\renewcommand{\S}{\mathbf{S}}

\begin{document}

\title{Constructions of sequential spaces}

\begin{abstract}
We introduce and study variable exponent $\ell^{p}$ spaces. These spaces will typically not be rearrangement-invariant 
but instead they enjoy a good local control of some required properties. Some interesting geometric examples are 
constructed by using these spaces.
\end{abstract}

\author[J. Talponen]{Jarno Talponen}
\address{University of Helsinki, Department of Mathematics and Statistics, Box 68, (Gustaf H\"{a}llstr\"{o}minkatu 2b) FI-00014 University 
of Helsinki, Finland}
\email{talponen@cc.helsinki.fi}

\maketitle

\section{Introduction}

In this paper we will introduce a natural schema for producing geometrically complicated Banach spaces with a $1$-unconditional basis.
The idea, however, is very simple, almost to the extent of being naive.
The resulting class of Banach spaces will be sufficiently flexible, so that within it one can construct easily Banach spaces 
with various combinations of suitable geometric and isomorphic properties.
 
Consider a map $p\colon \N\rightarrow [1,\infty]$. We will study variable exponent sequential spaces $\ell^{p(\cdot)}$ given formally by
\[\ell^{p(\cdot)}=\ldots(\ldots(((\R\oplus_{p(1)}\R)\oplus_{p(2)}\R)\oplus_{p(3)}\R)\oplus_{p(4)}\ldots\]
(This definition will be made rigorous shortly.) For example, for a constant function $p(\cdot)\equiv p\in [1,\infty]$ it holds that 
$\ell^{p(\cdot)}=\ell^{p}$ isometrically. These spaces differ from some other well-known variable exponent spaces in that 
the definition of the norm will be very explicit and the resulting sequential spaces will usually not be rearrangement-invariant.

By starting from $\ell^{p(\cdot)}$ type of space one can prove the following main result.
\begin{theorem}
There exists a Banach space with a $1$-unconditional Schauder basis that contains spaces $\ell^{p},\ 1\leq p<\infty,$
isomorphically, in fact almost isometrically. 
\end{theorem}

Recall the classical fact that $\ell^{p}$ and $\ell^{q}$ are non-isomorphic, when $p\neq q$. For this reason the above
result is perhaps surprising, as the claimed space is separable but must contain continuum many (mutually non-isomorphic) 
$\ell^{p}$-spaces. Recall that $C([0,1])$ is universal for separable spaces but it does not admit an unconditional basis, 
see \cite[p.24]{LTI}.

We will also study the basic properties of $\ell^{p(\cdot)}$ spaces. These spaces have some nice analogous properties
in comparison with the classical $\ell^{p}$ spaces. On the other hand, some unexpected problems arise as well, for example, 
we do not know in general whether $[(e_{n})]\subset \ell^{p(\cdot)}$ coincides with $c_{0}\cap \ell^{p(\cdot)}$.
As for the structural properties of $\ell^{p(\cdot)}$ spaces, it turns out for example that type and cotype become useful concepts 
and behave nicely in this setting.

\subsection{Preliminaries}
Given $t,s\in [0,\infty)$ and $p\in (0,\infty)$ we denote
\[t\boxplus_{p} s=(t^{p}+s^{p})^{\frac{1}{p}}\ \mathrm{and}\ t\boxplus_{\infty} s=\max(t,s).\]
Clearly $\boxplus_{p}$ gives a commutative semigroup for a fixed $p$. 
\begin{fact}\label{abcfact}
Let $0<p_{0}\leq p_{1}\leq \infty$ and $a,b,c\in [0,\infty)$. Then
\[(a\boxplus_{p_{0}}b)\boxplus_{p_{1}}c\leq a\boxplus_{p_{0}}(b\boxplus_{p_{1}}c).\]
\end{fact}
\begin{proof}
The claim is clearly equivalent to 
\begin{equation}\label{eq: abc}
((a^{p_{0}}+b^{p_{0}})^{\frac{p_{1}}{p_{0}}}+c^{p_{1}})^{\frac{p_{0}}{p_{1}}}\leq a^{p_{0}}+(b^{p_{1}}+c^{p_{1}})^{\frac{p_{0}}{p_{1}}},
\end{equation}
where $p_{0}\leq p_{1}<\infty$, and which holds as an equality for $a=0$. By substituting $u=a^{p_{0}}$, differentiating 
and using that $\frac{p_{0}}{p_{1}}-1\leq 0$ we obtain that
\begin{equation*}
\begin{array}{ll}
 & \frac{\partial}{\partial u}((u+b^{p_{0}})^{\frac{p_{1}}{p_{0}}}+c^{p_{1}})^{\frac{p_{0}}{p_{1}}}\leq 
\frac{\partial}{\partial u}((u+b^{p_{0}})^{\frac{p_{1}}{p_{0}}})^{\frac{p_{0}}{p_{1}}}\\
 & =1=\frac{\partial}{\partial u}(u+(b^{p_{1}}+c^{p_{1}})^{\frac{p_{0}}{p_{1}}})
\end{array}
\end{equation*}
for $u\geq 0$ and this yields \eqref{eq: abc}.
\end{proof}

Next we will give a precise definition for the variable-exponent $\ell^{p}$ spaces. 
Let $p\colon \N\rightarrow [1,\infty]$ be a map and $x\in \ell^{\infty}$.
We define semi-norms $|||\cdot|||_{k}$ on $\ell^{\infty}$ recursively by putting $|||x|||_{(1)}=|x_{1}| \boxplus_{p(1)}|x_{2}|$ and 
$|||x|||_{(k)}=|||x|||_{(k-1)}\boxplus_{p(k)} |x_{k+1}|$ for $k\in\N,\ k\geq 2$.   
Observe that $(|||x|||_{(k)})$ is a non-decreasing sequence for each $x\in \ell^{\infty}$. Hence 
we may put $\Phi\colon \ell^{\infty}\rightarrow [0,\infty],\ \Phi((x_{n}))=\lim_{k\rightarrow\infty}|||x|||_{(k)}$ for 
$x\in \ell^{\infty}$. Consider the vector space
\[\ell^{p(\cdot)}=\{(x_{n})\in \ell^{\infty}:\ \Phi((x_{n}))<\infty\},\]
which is equipped with the usual point-wise linear structure.
It is easy to see that the mapping $||\cdot||_{\ell^{p(\cdot)}}\stackrel{\cdot}{=}\Phi|_{\ell^{p(\cdot)}}$ is a norm on $\ell^{p(\cdot)}$.

Let $e_{1}=(1,0,0,0,\ldots);\ e_{2}=(0,1,0,0,\ldots);\ldots$ be the canonical unit vectors of $c_{0}$. We denote
$P_{k}\colon \ell^{p(\cdot)}\rightarrow [e_{1},\ldots, e_{k}],\ (x_{n})\mapsto (x_{1},\ldots, x_{k},0,0,0,\ldots)$. 
Observe that by the construction of the norm $||\cdot||_{\ell^{p(\cdot)}}$ it holds that 
$||P_{i}(x)||_{\ell^{p(\cdot)}}\leq ||P_{j}(x)||_{\ell^{p(\cdot)}}$, for $i,j\in \N,\ i\leq j,\ x\in \ell^{p(\cdot)}$.
Note that $||x||_{\ell^{p(\cdot)}}=\sup_{k\in\N}||P_{k}(x)||_{\ell^{p(\cdot)}}$ for $x\in \ell^{p(\cdot)}$ and that 
$P_{n}$ is a norm-$1$ projection for $n\in\N$. We will denote $Q_{n}\stackrel{\cdot}{=}\I-P_{n}$.

We will denote the Banach-Mazur distance of mutually isomorphic Banach spaces $\X$ and $\Y$ by 
\[d_{\mathrm{BM}}(\X,\Y)=\inf \{||T||\cdot ||T^{-1}||:\ T\colon \X\rightarrow \Y\ \mathrm{is\ an\ isomorphism}\}.\] 
Spaces $\X$ and $\Y$ are \emph{almost isometric} if $d_{\mathrm{BM}}(\X,\Y)=1$.
Recall that a Banach space $\X$ is \emph{contained almost isometrically} in a Banach space $\Z$ if
for each $\epsilon>0$ there is a subspace $\Y\subset\Z$ such that $d_{\mathrm{BM}}(\X,\Y)<1+\epsilon$.
Here we will often encounter sequential Banach spaces $\X$ and $\Y$ such that $\X$ contains $\Y$ almost isometrically in a 
such way that 
\begin{equation}\label{eq: ast}
\lim_{n\rightarrow \infty}d_{\mathrm{BM}}(Q_{n}(\X),\Y)=1.
\end{equation}

Recall that a Banach space $\X$ is an \emph{Asplund space} if any separable subspace of $\X$ has a separable dual.
Given a locally convex topology $\tau$ on $\X$, the space $\X$ is said to be \emph{$\tau$ locally uniformly rotund} 
($\tau$-LUR for short) if the following holds: For each sequence $(x_{n})\subset \S_{\X}$ such that 
$||x_{1}+x_{n}||\rightarrow 2$ as $n\rightarrow \infty$ it holds that $x_{n}\stackrel{\tau}{\longrightarrow} x_{1}$ as 
$n\rightarrow \infty$. If $\tau$ is the norm topology then we write LUR instead of $\tau$-LUR. The space $\X$ is 
\emph{midpoint locally uniformly rotund} (MLUR) if for each point $x\in \S_{\X}$ and sequences $(y_{n}),(z_{n})\subset \S_{\X}$
such that $\frac{1}{2}(y_{n}+z_{n})\rightarrow x$ it holds that $||y_{n}-z_{n}||\rightarrow 0$ as $n\rightarrow \infty$.
  
\section{Results}
The variable-exponent $\ell^{p}$ spaces can be used in constructing Banach spaces, which admit in some sense pathological, 
yet $1$-unconditional bases. Next we list examples of such results, the proofs of which are given subsequently. 

\begin{theorem}\label{mainthm}
The class of Banach spaces of the type $\ell^{p(\cdot)}$ contains a universal space up to almost isometric containment.
The analogous statement holds for spaces of the type $[(e_{n})]\subset \ell^{p(\cdot)}$.
\end{theorem}

\begin{theorem}\label{mainthm2}
Let $1<p<\infty$. Then there is a Banach space $\X$ with a $1$-unconditional basis $(f_{n})$ such that 
\begin{enumerate}
\item[(i)]{$\ell^{p}$ does not contain an isomorphic copy of $\X$.}
\item[(ii)]{For each strictly increasing subsequence $(i)\subset \N$ and $\epsilon>0$ there is a further subsequence $(i_{j})$
such that $[f_{i_{j}}:\ j\in\N]$ is $1+\epsilon$ isomorphic to $\ell^{p}$ in the sense of \eqref{eq: ast} 
via the equivalence of the bases.}
\end{enumerate}
\end{theorem}

\begin{theorem}\label{mainthm3}
Each space $\ell^{q(\cdot)}$ contains $\ell^{p}$ almost isometrically in the sense of \eqref{eq: ast} for some $p\in [1,\infty]$.
\end{theorem}

\subsection{Basic Properties}
It turns out that spaces of the type $\ell^{p(\cdot)}$ enjoy some basic properties similar to that of classical 
$\ell^{p}$ spaces.

\begin{proposition}\label{propo}
Let $p\colon \N\rightarrow [1,\infty]$. Then $\ell^{p(\cdot)}$ is a Banach space. 
Moreover, $(e_{n})_{n\in\N}$ is a $1$-unconditional basis of the space $[(e_{n})_{n\in\N}]$.
\end{proposition}
\begin{proof}
Clearly $[e_{1},\ldots ,e_{n}]\subset\ell^{p(\cdot)}$ is a Banach spaces for $n\in \N$. 
Let $(x^{(j)})_{j}\subset\ell^{p(\cdot)}$ be a Cauchy sequence. Note that $(x^{(j)})_{j}$ is bounded in $\ell^{\infty}$.
Since $P_{k}$ is a contractive projection for $k\in\N$, it holds that $(P_{n}(x^{(j)}))_{j}$ is a Cauchy sequence 
for a fixed $n$ and hence converges in $[e_{1},\ldots ,e_{n}]$. We conclude that $x^{(j)}\rightarrow x$ pointwise as $j\rightarrow\infty$,
for some $x\in \ell^{\infty}$. Let $\epsilon>0$. Since $(x^{(j)})_{k}\subset \ell^{p(\cdot)}$ is Cauchy
we get that there is $i_{0}\in \N$ such that $||x^{(j)}-x^{(i_{0})}||<\epsilon$ for $j\geq i_{0}$. In particular
\begin{equation}\label{eq: eps}
||P_{k}(x^{(j)}-x^{(i_{0})})||_{\ell^{p(\cdot)}}<\epsilon\quad \mathrm{for}\ n\in\N.
\end{equation}
By the definition of $\Phi$ and \eqref{eq: eps} we get that 
\[\Phi(x-x^{(i_{0})})=\sup_{k}||P_{k}(x-x^{(i_{0})})||_{\ell^{p(\cdot)}}<\epsilon.\]
This yields that $\Phi(x)\leq \Phi(x^{(i_{0})})+\epsilon<\infty$. Moreover, since $\epsilon$ was arbitrary
we get that $x^{(j)}\rightarrow x$ in $\ell^{p(\cdot)}$ as $j\rightarrow\infty$. This completes the 
claim that $\ell^{p(\cdot)}$ is complete.

To check the latter claim, let $x=(x_{n})_{n}\in [(e_{n})_{n}]$. We apply an auxiliary sequence 
$(y^{(j)})_{j}\subset \span((e_{n})_{n})$ such that $y^{(j)}\rightarrow x$ in $\ell^{p(\cdot)}$ as $j\rightarrow \infty$.
Set $k_{i}\in\N$ for $i\in\N$ such that $P_{k_{i}}(y^{(j)})=y^{(j)}$ for $j\leq i$. Observe that 
$||x-P_{k_{i}}(x)||_{\ell^{p(\cdot)}}\leq ||x-P_{k_{i}}(y^{(j)})||_{\ell^{p(\cdot)}}$ for $j\leq i$.
This yields that 
\[\sup_{i\geq j}||x-P_{k_{i}}(x)||_{\ell^{p(\cdot)}}\leq \sup_{i\geq j}||x-P_{k_{i}}(y^{(j)})||_{\ell^{p(\cdot)}}
=||x-y^{(j)}||_{\ell^{p(\cdot)}}\rightarrow 0\ \mathrm{as}\ j\rightarrow \infty.\]
We conclude that $(e_{n})_{n}$ is a Schauder basis. By the construction of $\Phi$ it holds that 
$||(x_{n})_{n}||_{\ell^{p(\cdot)}}=||(\theta(n)x_{n})_{n}||_{\ell^{p(\cdot)}}$ for any sequence of signs
$\theta\in \{-1,1\}^{\N}$. Thus $(e_{n})_{n}$ is a $1$-unconditional basis of $[(e_{n})_{n}]$.
\end{proof}

Clearly $\ell^{p(\cdot)}\cap c_{0}$ is a closed subspace of $\ell^{p(\cdot)}$ for any $p$, and hence it contains 
$[(e_{n})_{n\in\N}]$. We do not know if the space $\ell^{p(\cdot)}\cap c_{0}$ coincides with $[(e_{n})_{n\in\N}]\subset \ell^{p(\cdot)}$.
If $\limsup_{n\rightarrow \infty}p(n)<\infty$, then it is easy to check that $\ell^{p(\cdot)}\subset c_{0}$.

\begin{proposition}(H\"{o}lder inequality)
For a given sequence $p\colon \N\rightarrow [1,\infty]$ put $p^{\ast}\colon \N \rightarrow [1,\infty]$
such that $\frac{1}{p(n)}+\frac{1}{p^{\ast}(n)}=1$ for each $n\in\N$. If $(x_{n})$ and $(y_{n})$ are real valued 
sequences then   
\[\sum_{n\in\N}|x_{n}y_{n}|\leq ||(x_{n})||_{\ell^{p(\cdot)}}\ ||(y_{n})||_{\ell^{p^{\ast}(\cdot)}}.\]
\end{proposition}
\begin{proof}[Proof by induction.]
\end{proof}

\begin{proposition}
Let $p\colon \N\rightarrow [1,\infty]$ and let $\X=[(e_{n})]\subset \ell^{p(\cdot)}$.
Then $\X^{\ast}=\ell^{p^{\ast}(\cdot)}$, where the duality is given by $x^{\ast}(x)=\sum_{n} x_{n}^{\ast}(x_{n})$.
\end{proposition}
\begin{proof}
Recall that $(e_{n})$ is the $1$-unconditional basis of $[(e_{n})]$ and that $(\R\oplus_{p}\R)^{\ast}=\R\oplus_{p^{\ast}}\R$
for $1\leq p\leq\infty$. Note that $P_{m}^{\ast}(\X^{\ast})$ is isometric to $P_{m}(\ell^{p^{\ast}(\cdot)})$ for each $m\in \N$. 
For each $m$ we identify these spaces and the corresponding projection is denoted by 
$P_{[e_{1}^{\ast},\ldots,e_{m}^{\ast}]}^{\ast}\colon \X^{\ast}\rightarrow [e_{1}^{\ast},\ldots,e_{m}^{\ast}]$.
The projections of the type $P_{[e_{1}^{\ast},\ldots,e_{m}^{\ast}]}^{\ast}$ clearly commute. 
Since $(e_{n})$ is a Schauder basis of $[(e_{n})]$, we obtain that 
\begin{equation}\label{eq: duality}
x^{\ast}(x)=\lim_{m\rightarrow\infty} P_{[e_{1}^{\ast},\ldots,e_{m}^{\ast}]}^{\ast}(x^{\ast})(x)=\lim_{m\rightarrow\infty}x^{\ast}(P_{m}(x))
\end{equation}
for each $x\in \X$. 

There exists for each $x^{\ast}\in\X^{\ast}$ a unique vector $f_{x^{\ast}}\in \ell^{\infty}$ such that 
$f_{x^{\ast}}=\omega^{\ast}-\lim_{m\rightarrow \infty}P_{[e_{1}^{\ast},\ldots,e_{m}^{\ast}]}^{\ast}(x^{\ast})$, where we consider
$[e_{1}^{\ast},\ldots ,e_{m}^{\ast}]$ in the canonical way as a subspace of $\ell^{\infty}$.
Note that the continuity of given functionals $f,g\in \X^{\ast}$ yields that if $f(e_{n})=g(e_{n})$ for $n\in \N$, then $f-g=0$.
Thus, if $x^{\ast}\neq y^{\ast}$, then $f_{x^{\ast}}\neq f_{y^{\ast}}$. 

In fact $f_{x^{\ast}}\in \ell^{p^{\ast}(\cdot)}$ and moreover $||f_{x^{\ast}}||_{\ell^{p^{\ast}(\cdot)}}=||x^{\ast}||_{\X^{\ast}}$. 
Indeed, by using \eqref{eq: duality} and the basic properties of $(e_{n})$ we get that 
\begin{eqnarray*}
||x^{\ast}||_{\X^{\ast}}&=&\sup_{x\in \S_{\X}}|x^{\ast}(x)|=
\sup_{x\in \S_{\X}}\lim_{m\rightarrow\infty}|P_{[e_{1}^{\ast},\ldots,e_{m}^{\ast}]}^{\ast}(x^{\ast})(x)|\\
&=&\lim_{m\rightarrow\infty}||P_{[e_{1}^{\ast},\ldots,e_{m}^{\ast}]}^{\ast}(x^{\ast})||_{\ell^{p^{\ast}(\cdot)}}
=||f_{x^{\ast}}||_{\ell^{p^{\ast}(\cdot)}}.
\end{eqnarray*}
Hence $\X^{\ast}$ can be regarded as an isometric subspace of $\ell^{p^{\ast}(\cdot)}$ respecting the given duality.
Finally, the H\"{o}lder inequality gives that $\X^{\ast}=\ell^{p^{\ast}(\cdot)}$.
\end{proof}

\subsection{Almost isometric containment}
Next we will prove the results formulated previously. The arguments share some common auxiliary observations 
and we proceed by proving these facts. 

\begin{fact}
If $(n_{k})_{k}\subset \N$ is a sequence, then 
\[\{(x_{j})\in \ell^{p(\cdot)}|\ x_{j}\neq 0\implies j\in \{n_{1},n_{1}+1,n_{2}+1,n_{3}+1,\ldots\}\ \}\subset \ell^{p(\cdot)}\] 
is isometric to $\ell^{q(\cdot)}$, where $q(k)=p(n_{k})$ for $k\in\N$.
\end{fact}
This justifies the notation $\ell^{q(\cdot)}\subset \ell^{p(\cdot)}$. It is clear that if mappings 
$p_{1},p_{2}\colon \N\rightarrow [1,\infty]$ satisfy $p_{1}\leq p_{2}$ pointwise, 
then $\ell^{p_{2}(\cdot)}\subset \ell^{p_{1}(\cdot)}$.

For $p_{1},p_{2}\colon \N \rightarrow [1,\infty]$ and $k\in \N$ we denote by $p_{1}|^{(k)}p_{2}\colon \N\rightarrow [1,\infty]$, 
the map given by $p_{1}|^{(k)}p_{2}(n)=p_{2}(n)$ for $n\geq k$ and $p_{1}|^{(k)}p_{2}(n)=p_{1}(n)$ for $n<k$. 
Consider a finite sequence of mappings
\[\ell^{p_{1}|^{(k)}p_{2}}\ \stackrel{\psi_{k-1}}{\longrightarrow}\ \ell^{p_{1}|^{(k-1)}p_{2}}\ \stackrel{\psi_{k-2}}{\longrightarrow}\ 
\dots\ \stackrel{\psi_{k-i}}{\longrightarrow}\ \ell^{p_{1}|^{(k-i)}p_{2}}\ \stackrel{\psi_{k-i-1}}{\longrightarrow}\ \dots
\ \stackrel{\psi_{1}}{\longrightarrow}\ \ell^{p_{2}(\cdot)},\]
where $\psi_{j}=\I\colon \ell^{p_{1}|^{(j+1)}p_{2}}\rightarrow \ell^{p_{1}|^{(j)}p_{2}}$ for $1\leq j\leq k-1$.
Observe that 
\begin{equation}\label{eq: comment}
||\psi_{j}||=||\I\colon \ell^{p_{1}(j)}_{2}\rightarrow \ell^{p_{2}(j)}_{2}||\ \mathrm{and}\  
||\psi_{1}\circ\ldots\circ \psi_{k-1}||\leq \prod_{j=1}^{k-1}||\psi_{j}||.
\end{equation}

\begin{lemma}\label{lemma}
Let $p,q\colon\N\rightarrow [1,\infty]$ and $\epsilon>0$. 
If $\liminf_{n\rightarrow\infty}|q(k)-p(n)|=0$ for $k\in \N$, then there is a strictly increasing sequence 
$(n_{k})_{k}\subset\N$ such that $\phi\colon y_{k}\mapsto x_{n_{k}}$ defines an embedding
$\ell^{q(\cdot)}\rightarrow \ell^{p(\cdot)}$ such that $(1+\epsilon)^{-1}||y||\leq ||\phi(y)||\leq (1+\epsilon)||y||$
for $y\in \ell^{q(\cdot)}$. Moreover, $\phi$ is an embedding satisfying \eqref{eq: ast}.
\end{lemma}
\begin{proof}
Extract a subsequence $(n_{k})_{k}\subset\N$ such that 
$\prod_{k\in\N}||\I\colon \ell^{p(n_{k})}_{2}\rightarrow \ell^{q(k)}_{2}||\leq 1+\epsilon$ and
$\prod_{k\in\N}||\I\colon \ell^{q(k)}_{2}\rightarrow \ell^{p(n_{k})}_{2}||\leq 1+\epsilon$.
Define a linear map $\phi\colon \ell^{q(\cdot)}\rightarrow \ell^{\infty}$ by 
$\phi\colon a_{k+1}e_{k+1}\mapsto a_{k+1}e_{n_{k}+1}$.
Fix $y\in \ell^{q(\cdot)}$ and $m\in \N$. Next we will apply the preceding observations in \eqref{eq: comment}.  
We obtain that 
\[(1+\epsilon)^{-1}||P_{m}(y)||_{\ell^{q(\cdot)}}\leq 
||\phi(P_{m}(y))||_{\ell^{p(\cdot)}}\leq (1+\epsilon)||P_{m}(y)||_{\ell^{q(\cdot)}}.\] 
Thus by recalling the definition of the norms $||\cdot||_{\ell^{q(\cdot)}}$ and $||\cdot||_{\ell^{p(\cdot)}}$ we obtain that 
$\phi\colon \ell^{q(\cdot)}\rightarrow \ell^{p(\cdot)}$ is defined, and this is the claimed isomorphism.
By inspecting the construction of $\phi$ it is clear that also the latter part of the statement holds.
\end{proof}

\begin{proof}[Proof of Theorem \ref{mainthm}]
Enumerate $\Q\cap [1,\infty)=(q(n))_{n\in\N}$, where $q\colon \N\rightarrow \Q$. Let $\ell^{q(\cdot)}$ be the corresponding space. 
It is easy to see that, given $p\colon \N\rightarrow [1,\infty]$ and $k\in \N$, it holds that 
$\liminf_{n\rightarrow \infty}|q(n)-p(k)|=0$. Lemma \ref{lemma} yields that $\ell^{q(\cdot)}$ contains $\ell^{p(\cdot)}$ 
almost isometrically, so that the first claim holds.

In the latter claim $[(e_{n})_{n}]\subset \ell^{q(\cdot)}$ is a suitable universal space. 
Indeed, according to Proposition \ref{propo} $(e_{n})_{n}$ is an unconditional basis
and the $1+\epsilon$-isomorphism $\phi$ appearing in the proof of Lemma \ref{lemma} takes $(e_{n})\subset \ell^{p(\cdot)}$
to $(f_{n_{1}},f_{n_{1}+1}, f_{n_{2}+1},\ldots)\subset \ell^{q(\cdot)}$. 
\end{proof}

In the above proof the space $\ell^{q(\cdot)}$ contains $\ell^{\infty}$ almost isometrically in the sense of \eqref{eq: ast}. 
Hence it is easy to see that $\ell^{q(\cdot)}/c_{0}$ contains $\ell^{\infty}/ c_{0}$ isometrically.
It is a classical fact that $\ell^{\infty}$ contains separable Banach spaces isometrically, and using the same argument,
so does $\ell^{\infty}/c_{0}$. We conclude that $\ell^{q(\cdot)}/ c_{0}$ contains separable spaces isometrically.

\begin{proof}[Proof of Theorem \ref{mainthm2}]
Let $1<p<\infty$. Let $(r_{i})\subset (1,\infty)$ be a sequence such that $r_{i}\rightarrow p$ as $i\rightarrow \infty$, 
$r_{n}<p$ for $n\in\N$ if $p\leq 2$ and $r_{n}>p$ for $n\in\N$ if $p>2$. 
The space $\ell^{p}$ has type $p$ if $p\leq 2$ and cotype $p$ if $p\geq 2$, and in both cases, given $i\in \N$, $\ell^{p}$ 
does not contain $\ell^{r_{i}}_{n}$s uniformly (see e.g. \cite{Pi}).
Hence we may pick for each $i\in\N$ a number $j_{i}\in \N$ such that 
\[\inf_{E}d_{\mathrm{BM}}(E,\ell^{r_{i}}_{j_{i}})>i,\]
where the infimum is taken over $j_{i}$-dimensional subspaces $E$ of $\ell^{p}$.
Define $q\colon \N\rightarrow [1,\infty]$ by putting $q(n)=r_{1}$ for $1\leq n\leq r_{1}$ and 
$q(n)=r_{l}$ for $\sum_{i<l}r_{i}<n\leq \sum_{i\leq l}r_{i}$ and $l>1$. Then it follows from the selection of 
$q\colon \N\rightarrow [1,\infty]$ that $\ell^{q(\cdot)}$ does not embed linearly into $\ell^{p}$.

Consider the canonical unit vectors $(f_{n})$ of $\ell^{q(\cdot)}$. According to Proposition \ref{propo} 
$(f_{n})$ is an unconditional basis of $\X=[(f_{n})]$. Since $q(n)\rightarrow p$ as $n\rightarrow\infty$, 
an application of Lemma \ref{lemma} yields the claim.
\end{proof}

\begin{proof}[Proof of Theorem \ref{mainthm3}]
Let $q\colon \N\rightarrow [1,\infty]$. Since $[1,\infty]$ is a compact metrizable space, therefore sequentially compact, 
there exists a subsequence $(q(n_{k}))_{k\in\N}\subset (q(n))_{n\in\N}$ convergent in $[1,\infty]$.
Let $p=\lim_{k\rightarrow \infty}q(n_{k})\in [1,\infty]$. Thus we may apply Lemma \ref{lemma} to conclude that 
$\ell^{q(\cdot)}$ contains $\ell^{p}$ almost isometrically.
\end{proof}

\begin{theorem}
Let $p\colon \N\rightarrow [1,\infty]$. Then the following conditions are equivalent:
\begin{enumerate}
\item[(1)]{$\ell^{p(\cdot)}$ is reflexive.}
\item[(2)]{$\ell^{p(\cdot)}$ is superreflexive.}
\item[(3)]{$\liminf_{n\rightarrow\infty}p(n)>1$ and $\limsup_{n\rightarrow \infty}p(n)<\infty$.}
\item[(4)]{$\ell^{p(\cdot)}$ is an Asplund space.}
\end{enumerate}
\end{theorem}

Let us observe before passing to the proof that $\ell{1}$ has the RNP but is not a reflexive space and thus 
one cannot replace Asplund by RNP in $(4)$.
\begin{proof}
Recall the well-known characterization of superreflexive spaces, due to Enflo, that a space is superreflexive if and only if it
is isomorphic to a space both uniformly convex and uniformly smooth. Clearly (2)$\implies$(1). 
By using Lemma \ref{lemma} we obtain that if
$\liminf_{n\rightarrow\infty}p(n)=1$ (resp. $\limsup_{n\rightarrow \infty}p(n)=\infty$), then 
$\ell^{p(\cdot)}$ contains $\ell^{1}$ (resp. $\ell^{\infty}$) almost isometrically, and thus fails to be
reflexive. Hence (1)$\implies$(3) holds. Similarly (4)$\implies$(3) holds and it is well-known that 
reflexive spaces are Asplund.

It suffices to verify direction (3)$\implies$(2), so let us assume that (3) holds. 
Then there exists $p_{0}\in (1,\liminf_{n\rightarrow \infty}p(n)]$, 
$q_{0}\in [\limsup_{n\rightarrow\infty}p(n),\infty)$ and $k_{0}\in \N$ such that 
$p(n)\in [p_{0},q_{0}]$ for $n\geq k_{0}$. Define $\widetilde{p}\colon \N\rightarrow [p_{0},q_{0}]$
by $\widetilde{p}(n)=\min(q_{0},\max(p(n),p_{0}))$ for $n\in\N$.
Note that $\span(e_{1},\ldots,e_{k})$ is a bicontractively complemented subspace regardless of 
whether considered as being contained in $\ell^{p(\cdot)}$ or $\ell^{\widetilde{p}(\cdot)}$. It follows that 
the identity mapping $I\colon \ell^{p(\cdot)}\rightarrow \ell^{\widetilde{p}(\cdot)}$ is an isomorphism.
Thus our task is reduced to proving that $\ell^{\widetilde{p}(\cdot)}$ is superreflexive.

We will require the notions of upper $p$-estimate and lower $q$-estimate of Banach lattices.
If $\X$ is a Banach lattice and $1\leq p\leq q<\infty$ then the upper $p$-estimate and the lower $q$-estimate, respectively, 
are defined (for the relevant multiplicative constants being $1$) as follows:\\
\begin{equation*}
\begin{array}{l}
\left|\left|\sum_{1\leq i\leq n} x_{i}\right|\right|\leq (\sum_{1\leq i\leq n}||x_{i}||^{p})^{\frac{1}{p}},\phantom{\bigg |}\\
\left|\left|\sum_{1\leq i\leq n} x_{i}\right|\right|\geq (\sum_{1\leq i\leq n}||x_{i}||^{q})^{\frac{1}{q}},\phantom{\bigg |}
\end{array}
\end{equation*}
respectively, for any disjoint vectors $x_{1},\ldots, x_{n}\in \X$.
We will apply the fact that a Banach lattice, which satisfies an upper $p$-estimate and a lower $q$-estimate 
for some $1<p<q<\infty$ is isomorphic to a Banach space both uniformly convex and uniformly smooth, see \cite[1.f.1, 1.f.7]{LTII}.
Note that in such case $\X$ is superreflexive and that the space $\ell^{\widetilde{p}(\cdot)}$, having a $1$-unconditional basis, 
carries a natural Banach lattice structure. Thus, it suffices to show that $\ell^{\widetilde{p}(\cdot)}$ satisfies
upper and lower estimates for $p_{0}$ and $q_{0}$, respectively. 

Denote by $P_{m}\colon \ell^{\widetilde{p}(\cdot)}\rightarrow \span(e_{1},\ldots,e_{m})$ the natural projection
preserving the first $m$ coordinates. To check that $\ell^{\widetilde{p}(\cdot)}$ satisfies the upper $p_{0}$-estimate, 
let $x_{1},\ldots,x_{n}\in \ell^{\tilde{p}(\cdot)}$ be disjoint vectors. We claim that
\begin{equation*}
\begin{array}{l}
\left|\left|\sum_{1\leq i\leq n} x_{i}\right|\right|_{\ell^{\widetilde{p}(\cdot)}}\leq \left(\sum_{1\leq i\leq n}||x_{i}||_{\ell^{\widetilde{p}(\cdot)}}^{p_{0}}\right)^{\frac{1}{p_{0}}}.
\end{array}
\end{equation*}
Indeed, assume to the contrary that this does not hold and let $m\in \N$ be the least natural number such that
\begin{equation*}
\begin{array}{l}
\left|\left|\sum_{1\leq i\leq n} P_{m}x_{i}\right|\right|_{\ell^{\widetilde{p}(\cdot)}}> \left(\sum_{1\leq i\leq n}||P_{m}x_{i}||_{\ell^{\widetilde{p}(\cdot)}}^{p_{0}}\right)^{\frac{1}{p_{0}}}.
\end{array}
\end{equation*}
Clearly $m>1$. We may assume without loss of generality that $(P_{m}-P_{m-1})x_{n}\stackrel{\cdot}{=}y_{m}\neq 0$. 
It follows from the disjointness of $x_{1},\ldots,x_{n}$ that $(P_{m}-P_{m-1})x_{i}=0$ for $1\leq i<n$.
Observe that 
\begin{equation*}
\begin{array}{l}
\left|\left|\sum_{1\leq i\leq n} P_{m-1}x_{i}\right|\right|_{\ell^{\widetilde{p}(\cdot)}}\leq \left(\sum_{1\leq i\leq n}||P_{m-1}x_{i}||_{\ell^{\widetilde{p}(\cdot)}}^{p_{0}}\right)^{\frac{1}{p_{0}}}
\end{array}
\end{equation*}
by the selection of $m$. Then 
\begin{equation*}
\begin{array}{lll}
 & & ||\sum_{1\leq i\leq n} P_{m}x_{i}||_{\ell^{\widetilde{p}(\cdot)}}=||\sum_{1\leq i\leq n} P_{m-1}x_{i}||_{\ell^{\widetilde{p}(\cdot)}}\boxplus_{\tilde{p}(m-1)}|y_{m}|\phantom{\bigg |}\\
 & \leq & (((\sum_{1\leq i<n}||P_{m}x_{i}||_{\ell^{\widetilde{p}(\cdot)}}^{p_{0}})^{\frac{1}{p_{0}}}\boxplus_{p_{0}} ||P_{m-1}x_{n}||)^{p_{0}})^{\frac{1}{p_{0}}}\boxplus_{\tilde{p}(m-1)}|y_{m}|\phantom{\bigg |}\\
 & \leq & (((\sum_{1\leq i<n}||P_{m}x_{i}||_{\ell^{\widetilde{p}(\cdot)}}^{p_{0}})^{\frac{1}{p_{0}}}\boxplus_{p_{0}}(||P_{m-1}x_{n}||\boxplus_{\tilde{p}(m-1)}|y_{m}|)\phantom{\bigg |}\\ 
 & = & (\sum_{1\leq i\leq n}||P_{m}x_{i}||_{\ell^{\widetilde{p}(\cdot)}}^{p_{0}})^{\frac{1}{p_{0}}}.\phantom{\bigg |}
\end{array}
\end{equation*}
Above we applied the selection of $m$ and Fact \ref{abcfact}.
Thus we arrive in a contradiction, which means that $\ell^{\widetilde{p}(\cdot)}$ satisfies the upper $p_{0}$-estimate.
The lower $q_{0}$-estimate is checked in a similar manner.
\end{proof}

\section{Final remarks}

A question was raised in \cite[p.174]{KR} whether each MLUR Banach space $\X$ is LUR. It has been established by now 
that this is not the case (see e.g. \cite{Haydon}). Next we will give a rather simple and natural example, which is 
related to these convexity conditions. 

\begin{proposition}\label{ex1}
Let $p\colon \N\rightarrow (1,\infty)$ be such that $\prod_{k\in \N}||\I\colon \ell^{p(k)}_{2}\rightarrow \ell^{1}_{2}||<2$.
Then $\ell^{p(\cdot)}$ satisfies the following condition: Given $x\in \S_{\ell^{p(\cdot)}}$ and sequences 
$(y_{n}),(z_{n})\subset \B_{\ell^{p(\cdot)}}$ such that $\frac{1}{2}(y_{n}+z_{n})\rightarrow x$
and $||P_{k}y_{n}||,||P_{k}z_{n}||\rightarrow ||P_{k}x||$ as $n\rightarrow \infty$ for each $k\in \N$, 
then $||y_{n}-z_{n}|| \rightarrow 0$ as $n\rightarrow \infty$.
However, $\ell^{p(\cdot)}$ is not  $\omega$-LUR. 
\end{proposition}

\begin{proof}
Clearly the sequence of canonical unit basis vectors $(e_{n})\subset \ell^{1}\subset \ell^{p(\cdot)}$ 
is a Schauder basis for $\ell^{p(\cdot)}$, since $\ell^{p(\cdot)}$ and $\ell^{1}$ are isomorphic by the construction of $p(\cdot)$. 
We denote by $P_{k}\colon \ell^{p(\cdot)}\rightarrow [e_{1},\ldots,e_{k}]$ the projection given by 
$\sum_{i\in\N}a_{i}e_{i}\longmapsto \sum_{i=1}^{k}a_{i}e_{i}$ for $k\in\N$. 

Fix $x\in\S_{\ell^{p(\cdot)}}$ and $(y_{n}),(z_{n})\subset\B_{\ell^{p(\cdot)}}$ as in the assumptions and 
we aim to show that $y_{n}-z_{n}\rightarrow 0$ as $n\rightarrow\infty$. 
Since $(y_{n})$ and $(z_{n})$ are otherwise arbitrary, it suffices to check without loss of generality that
there exists a subsequence $(n_{k})\subset \N$ such that $y_{n_{k}}-z_{n_{k}}\rightarrow 0$ as $k\rightarrow\infty$.

By the continuity of $P_{k}$ we obtain that $\frac{1}{2}P_{k}(z_{n}+y_{n})\rightarrow P_{k}(x)$ for $k\in\N$.
Observe that $[e_{1},\ldots,e_{k}]\subset \ell^{p(\cdot)}$ is a uniformly convex subspace for each $k\in\N$.
Thus $P_{k}(z_{n}),P_{k}(y_{n})\rightarrow P_{k}(x)$ as $n\rightarrow \infty$ for all $k\in\N$.
Hence one can pick a sequence $(n_{k})\subset \N$ such that $P_{k}(y_{n_{k}})-P_{k}(z_{n_{k}})\rightarrow 0$
as $k\rightarrow \infty$. 

Similarly as earlier in \eqref{eq: comment}, one can consider $\I\colon \ell^{p(\cdot)}\rightarrow \ell^{1}$ formally as 
$\psi_{1}\circ \psi_{2}\circ\ldots\colon \ell^{p(\cdot)}\rightarrow \ell^{1}$. By applying the fact that 
$\lim_{i\rightarrow\infty}\lim_{k\rightarrow\infty}||\psi_{k-i}\circ \ldots\circ \psi_{k}||=1$
we obtain that the sequence of mappings $R_{k}\colon \ell^{p(\cdot)}\rightarrow \R$, given by
$R_{k}(x)=||P_{k}(x)||+||(\I-P_{k})(x)||$ for $k\in\N$, satisfies that $||R_{k}||\rightarrow 1$ as $k\rightarrow\infty$.

This is applied as follows. Since $||P_{k}(y_{n_{k}})||,||P_{k}(z_{n_{k}})||\rightarrow 1$ as $k\rightarrow\infty$, we obtain that
$||(\I-P_{k})(y_{n_{k}}-z_{n_{k}})||\rightarrow 0$ as $k\rightarrow\infty$. 
Hence we obtain that 
\[||z_{n_{k}}-y_{n_{k}}||_{\ell^{p(\cdot)}}\leq (||P_{k}(z_{n_{k}}-y_{n_{k}})||_{\ell^{p(\cdot)}}+||(\I-P_{k})(z_{n_{k}}-y_{n_{k}})||_{\ell^{p(\cdot)}})\stackrel{k\rightarrow\infty}{\longrightarrow} 0.\] 
Consequently $\ell^{p(\cdot)}$ satisfies the first claim.

One can pick a strictly increasing sequence $(n_{i})\subset \N$ such that 
$||(1-2^{-i},1-2^{-i})||_{\ell^{p^{\ast}(n_{i})}_{2}}\leq 1-2^{-i-1}$ for all $i\in\N$. Define a sequence
$(x_{n})$ by $x_{n_{i}+1}=1-2^{-i}$ for all $n_{i}$ and $x_{n}=0$ for all $n\in\N\setminus (n_{i}+1)_{i}$. 
Observe that $(x_{n})\in\B_{\ell^{p^{\ast}(\cdot)}}$ and that $x_{1}=0$. According to the H\"{o}lder inequality
$f\colon \ell^{p(\cdot)}\rightarrow \R,\ (y_{n})\mapsto \sum_{n\in\N}x_{n}y_{n}$
is defined and $f\in \B_{(\ell^{p(\cdot)})^{\ast}}$. 

Finally, observe that $||e_{1}+e_{n_{i}+1}||_{\ell^{p(\cdot)}}\rightarrow 2$ as $i\rightarrow\infty$. 
However, $\lim_{i\rightarrow\infty}f(e_{n_{i}+1})=1\neq 0=f(e_{1})$. This means that $\ell^{p(\cdot)}$ is not $\omega$-LUR.
\end{proof}

Observe that for each $\epsilon>0$ the space $\ell^{p(\cdot)}$ above can be additionally defined so that
its Banach-Mazur distance to $\ell^{1}$ is less than $1+\epsilon$. In fact,
\[\lim_{n\rightarrow \infty}d_{\mathrm{BM}}(Q_{n}(\ell^{p(\cdot)}),\ell^{1})=1.\]

Finally we reiterate the open problems that have arised in this note: \\
$\bullet$ We do not know if the space in Example \ref{ex1} is MLUR.\\
$\bullet$ Does $\ell^{p(\cdot)}\cap c_{0}$ always coincide with $[(e_{n})]\subset \ell^{p(\cdot)}$?\\
$\bullet$ Given $p\colon \N\rightarrow (0,\infty)$, is $\ell^{p(\cdot)}$ necessarily strictly convex or smooth?

\end{document}